\documentclass[11pt,reqno,a4paper]{amsart}

\usepackage{etoolbox}
\usepackage{amsmath}
\usepackage{enumerate}
\usepackage{amssymb}
\usepackage[initials]{amsrefs}
\usepackage[all]{xy}
\usepackage{graphicx,latexsym,amsbsy,amscd,amsthm,epsfig,verbatim}
\usepackage{amsfonts}
\usepackage{comment}
\usepackage{courier} 
\usepackage[T1]{fontenc}                                                    
\usepackage{hyperref}
\usepackage{slashed}
\usepackage{tikz-cd}

\numberwithin{equation}{section}

\newcommand\subsetsim{\mathrel{%
\ooalign{\raise0.2ex\hbox{$\subset$}\cr\hidewidth\raise-0.8ex\hbox{\scalebox{0.9}{$\sim$}}\hidewidth\cr}}}

\newcommand{\signum}{\operatorname{sign}}

\newtheorem{theorem}{Theorem}[section]

\newtheorem{proposition}[theorem]{Proposition}

\newtheorem{claim}[theorem]{Claim}
\theoremstyle{definition}
\newtheorem{definition}[theorem]{Definition}

\newtheorem{remark}[theorem]{Remark}

\patchcmd{\subsection}{-.5em}{.5em}{}{}
\patchcmd{\subsubsection}{-.5em}{.5em}{}{}

\title[Surface bundles over surfaces: new inequalities]{Surface bundles over surfaces: new inequalities between signature, simplicial volume and Euler characteristic}

\author{Michelle Bucher}
\email{michelle.bucher-karlsson@unige.ch}
\address{Section de Math\'ematiques, Universit\'e de Gen\`eve}

\author{Caterina Campagnolo}
\email{caterina.campagnolo@ens-lyon.fr}
\address{Unité de mathématiques pures et appliquées, ENS Lyon}

\subjclass[2010]{Primary 55R10; Secondary 57R22}

\keywords{Surface bundle, simplicial volume, signature}

\begin{document}
\begin{abstract}
We present three new inequalities tying the signature, the simplicial volume and the Euler characteristic of surface bundles over surfaces. Two of them are true for any surface bundle, while the third holds on a specific family of surface bundles, namely the ones that arise through ramified coverings. These are among the main known examples of bundles with non-zero signature. 
\end{abstract}

\maketitle

\section{Introduction}
Surface bundles over surfaces form an interesting family of $4$-manifolds that give rise to several questions: for example, do such manifolds with non-zero signature exist? If yes, which values does the signature take? What are the minimal base and fibre genera required to achieve a given signature? The relations and inequalities between signature and Euler characteristic of surface bundles have been widely studied, notably by Bryan, Catanese, Donagi, Endo, Korkmaz, Kotschick, Ozbagci, Rollenske, Stipsicz \cites{bryan-donagi, bryan-donagi-stipsicz, catanese-rollenske, endo, endo-korkmaz-kotschick-ozbagci-stipsicz, kotschick}.

In the present note we add the comparison to the simplicial volume of the total space, using tools from bounded cohomology. The simplicial volume can act as a bridge between the two other invariants, signature and Euler characteristic (for the definition of simplicial volume, see Subsection \ref{def vol simp}).

For any surface bundle $E$ over a surface, the best known inequality between the signature $\sigma(E)$ and the Euler characteristic $\chi(E)$ is due to Kotschick \cite{kotschick}:
$$2|\sigma(E)|\leq \chi(E).$$
Kotschick also obtained the stronger inequality $3|\sigma(E)|<\chi(E)$ in some special cases \cite{kotschick-3}.
The first author's work on simplicial volume of surface bundles \cite{bucher} produced an inequality between simplicial volume and Euler characteristic of aspherical surface bundles:
$$6\chi(E)\leq \|E\|.$$
We compare here the signature to the simplicial volume of general surface bundles over surfaces and obtain:

\begin{theorem}\label{thm:sign vol simp}
Let $E$ be an oriented surface bundle over a surface, with closed oriented base and fibre. Then
$$36|\sigma(E)|\leq \|E\|.$$
\end{theorem}
Observe that this is stronger than the combination of Kotschick's and the first author's inequalities, which only give $12|\sigma(E)|\leq \|E\|$, or $18|\sigma(E)|< \|E\|$ in the special cases of \cite{kotschick-3}. The inequality of Theorem \ref{thm:sign vol simp} is also strictly stronger than the value produced by the up to now best example \cite{catanese-rollenske}*{Theorem A}, which is $27|\sigma(E)|\leq \|E\|$.

The simplicial volume remains very hard to compute explicitly. In fact, the exact values in non-vanishing cases are known only for hyperbolic manifolds (due to Gromov-Thurston \cites{gromov, thurston}) and for locally $(\mathbb{H}^2\times\mathbb{H}^2)$-manifolds, so in particular for products of surfaces \cite{bucherH^2xH^2}.

We can give a lower bound on $\|E\|$ under the form of the $\ell_1$-norm of a distinguished $2$-homology class: 

\begin{proposition}\label{prop: N et E}
Let $E$ be an oriented surface bundle over a surface, with closed oriented base and fibre.
Let $\left[N\right]$ be the Poincar\'e dual of the Euler class of the tangent bundle along the fibre of $E$
. Then
$$\|\left[N\right]\|_1\leq \frac{1}{3}\|E\|.$$
\end{proposition}
The tangent bundle along the fibre will be defined in Subsection \ref{def classe euler}.
Observe that the dual of this Euler class can be represented by a subsurface of $E$, hence once we know its minimal genus we will be able to compute its $\ell_1$-norm. Unfortunately for now the known lower bounds on $\|\left[N\right]\|_1$ do not produce better inequalities for $\|E\|$ than the already existing ones.

Signatures remain, analogously to simplicial volume, quite hard to calculate for general surface bundles and are essentially only computed for bundles coming from specific constructions: differences of Lefschetz fibrations or ramified coverings. More recently, Baykur used yet another method in \cite{baykur}, namely horizontal and vertical stabilizations, and obtained infinite families of surface bundles with non-zero signature. We will specialise to the examples arising through ramified coverings (see Section \ref{ramifie} for the definition and notations) and prove:

\begin{theorem}\label{thm:fibre ramifie}
Let $E$ be a surface bundle as in Section \ref{ramifie}. Then
$$\|E\|\geq 6\chi(E)+6|\chi(\Sigma')|(d-1),$$
where $\Sigma'$ is the base of the bundle 
and $d$ is the degree of the ramified covering.
\end{theorem}
Remark that this improves the inequality $\|E\|\geq 6\chi(E)$ of the first author. It constitutes the first example of surface bundles over surfaces for which the strict inequality $\|E\|> 6\chi(E)$ is shown.

In the next section we recall the definitions of the invariants under consideration and the main tools to compute them. We devote Section \ref{preuve sign vol simp} to the proof of Theorem \ref{thm:sign vol simp} and Section \ref{preuve N et E} to the proof of Proposition \ref{prop: N et E}. The bundles related to ramified coverings will be treated in Section \ref{ramifie}.

\vspace{0.1cm}
\emph{Acknowledgements.}
This research was supported by the Swiss National Science Foundation. The second author is grateful to the first author, her doctoral advisor, for introducing her to the beautiful topics and techniques of bounded cohomology and simplicial volume. The authors thank Inan\c{c} Baykur for pointing out his interesting construction of surface bundles with non-zero signature. They also thank Dieter Kotschick for noticing an inaccuracy in a previous statement of Proposition \ref{prop: N et E} and Pierre de la Harpe for useful comments on an earlier version. The Karlsruhe Institute of Technology, former home institution of the second author, kindly takes care of the Open Access publication fee.

\section{Definition of the invariants}\label{def back}
In what follows we study oriented surface bundles over surfaces $F\hookrightarrow E \stackrel{\pi}{\twoheadrightarrow} B$, where both $F$ and $B$, and hence $E$, are closed.

While the \emph{Euler characteristic} does not need to be redefined, let us just recall that it is multiplicative in the base and the fibre of a bundle, that is it satisfies
$$\chi(E)=\chi(F)\chi(B).$$
In particular all the bundles with same base and fibre have the same Euler characteristic.
\subsection{Signature}
The \emph{signature} of a closed connected oriented $4k$-manifold $M$, where $k\in \mathbb{N}$, is defined as follows.

Consider the bilinear form induced by the cup product on the middle-dimensional cohomology groups:
$$\begin{array}{rrcl}
\cup\colon & H^{2k}(M, \mathbb{Z})\times H^{2k}(M, \mathbb{Z})& \longrightarrow & H^{4k}(M, \mathbb{Z})\cong \mathbb{Z}\\
& (\alpha, \beta)& \longmapsto & \alpha\cup \beta.
\end{array}$$
As $\alpha\cup\beta=(-1)^{2k\cdot 2k}\beta\cup\alpha=\beta\cup\alpha$, the form is symmetric. Thus all its eigenvalues are real, and we can compute its signature in $\mathbb{Z}$ as the number of positive eigenvalues $b_2^+(M)$ minus the number of negative eigenvalues $b_2^-(M)$. The $0$ eigenvalues are neglected.

The signature of $M$, denoted by $\sigma(M)$, is the signature of the above bilinear form.

\subsection{Simplicial volume}\label{def vol simp}
Let $X$ be a topological space.

One can define a semi-norm on homology classes in the singular homology: let $\zeta\in H_k(X, \mathbb{R})$ be a homology class. Then 
$$\|\zeta\|_1=\mathrm{inf}\left\{\sum_i |a_i| \, \vline\, \zeta=\left[\sum_ia_i\sigma_i\right]\in H_k(X, \mathbb{R})\right\},$$
where $\sigma_i\colon \Delta^k\rightarrow X$ denotes a singular simplex of dimension $k$.
We call this semi-norm the \emph{$\ell_1$-norm}.

The \emph{simplicial volume} of a closed oriented manifold $M$ of dimension $n$ is then defined as the $\ell_1$-norm of its (real) fundamental class $[M]$, 
$$\|M\|=\inf\left\{\sum_i |a_i|\,\vline\, \left[M\right]=\left[\sum_i a_i\sigma_i\right]\in H_n(M, \mathbb{R})\right\}.$$

This invariant was introduced by Gromov in \cite{gromov}. The simplicial volume has many facets: among others, it is a topological measure of the complexity of a manifold, it gives restrictions on the geometries a manifold can carry and admits immediate degree theorems.

We will also need the norm commonly used in the theory of bounded cohomology, but which we consider on standard singular cohomology classes. Let $\beta\in H^k(X, \mathbb{R})$ be a cohomology class. The \emph{norm} of $\beta$ (indeed a semi-norm) is defined as the infimum of the sup norm of all cochains representing $\beta$:
$$\|\beta\|=\mathrm{inf}\left\{\|b\|_\infty \,\vline\, b\in C^k(X, \mathbb{R}), [b]=\beta\right\}.$$
Note that it is possible that $\|b\|_\infty=\infty$ for every such $b$ and in particular that $\|\beta\|=\infty$.

We will use the following relationship between $\ell_1$-norm and sup norm:
\begin{proposition}[\cite{benedetti-petronio}, Proposition F.2.2]\label{norme-vol simp}
Let $\beta\in H^k(X,\mathbb{R}), \zeta\in H_k(X,\mathbb{R})$ as above. Then 
$$\frac{|\left\langle\beta, \zeta\right\rangle|}{\|\beta\|}\leq \|\zeta\|_1.$$
If $M$ is an oriented compact $n$-dimensional manifold and $\beta\in H^n(M, \mathbb{R})$ is a cohomology class of degree $n$, then 
$$\frac{|\left\langle \beta, \left[M\right]\right\rangle|}{\|\beta\|}=\|M\|.$$
\end{proposition}
\subsection{The Euler class}\label{def classe euler}
Let $E$, as above, be a surface bundle $F\hookrightarrow E \stackrel{\pi}{\twoheadrightarrow} B$. One defines its tangent bundle along the fibre as
$$T\pi=\{v\in TE\,\mid\, \pi_*(v)=0\}.$$
As an oriented vector bundle, it has an Euler class. We call it the \emph{Euler class} of the bundle $E$ and denote it by $e\in H^2(E, \mathbb{Z})$ --- not to be confused with the Euler class in top degree of $E$. Its Poincar\'e dual $e\cap [E]\in H_2(E, \mathbb{Z})$ will be denoted by $[N]$. Note that $[N]$, as a degree $2$ homology class in a $4$-manifold, is representable by a subsurface of $E$ (see for example \cite{hopf}).

The class $e$ has a quite explicit representative, which can be described as follows (see also \cite{bucher}*{Section 3}).
The holonomy morphism of the bundle $F\hookrightarrow E\stackrel{\pi}{\twoheadrightarrow} B$ gives rise to the following diagram:
$$
\xymatrix
{
\pi_1(E)\ar[r]\ar[d] \ar@/^1.5pc/[rr]^\phi& \mathcal{M}_{g, *}\ar[r]^-\rho\ar[d]& \mathrm{Homeo}_+(S^1)\\
\pi_1(B) \ar[r] & \mathcal{M}_g
}
$$
where $g$ is the genus of the fibre $F$ and $\mathcal{M}_g$ its mapping class group, while $\mathcal{M}_{g, *}$ denotes the group of mapping classes of $F$ fixing a given base point (see \cite{morita-article}*{Paragraphs 2 and 4}).

In $H^2(\mathrm{BHomeo}_+(S^1), \mathbb{Z})\cong H^2(\mathrm{Homeo}_+(S^1), \mathbb{Z})$ we have the Euler class $\chi$ that classifies flat $S^1$-bundles. Passing through the isomorphism $H^2(E, \mathbb{Z})\cong H^2(\pi_1(E), \mathbb{Z})$, true for aspherical $E$, \cite{morita-article}*{Proposition 4.1} gives us $e=\phi^*(\chi)$.

Now the class $\chi$ can be represented by half the \emph{orientation cocycle} on the circle \cite{morita-article}*{Proposition 4.3}. 
The orientation cocycle $Or$ is a $2$-cocycle defined as follows. Choose a point $x\in S^1$. Then:
$$
Or \, : \, \left\{
\aligned
(\mathrm{Homeo}_+(S^1))^3 \hskip.3cm &\longrightarrow \hskip.4cm \mathbb Z
\\
(g_0,g_1,g_2) \hskip.2cm &\longmapsto \hskip.2cm
\left\{
\aligned
&\phantom{-}1 \hskip.2cm \text{if $g_0x,g_1x,g_2x$ are distinct and positively oriented},
\\
&\phantom{-}0 \hskip.2cm \text{if two points among $g_0x,g_1x,g_2x$ coincide},
\\
&{-}1 \hskip.2cm \text{if $g_0x,g_1x,g_2x$ are distinct and negatively oriented}.
\endaligned
\right.
\endaligned
\right.
$$
It is alternating and its norm as a cocycle is obviously $1$. Moreover the cohomology class it defines does not depend on the choice of $x$.
							
Therefore the class $e$ can be represented by $\frac{1}{2}\phi^*(Or)$, so that it has an alternating representative and has norm $\|e\|\leq\frac{1}{2}$.

The signature of a surface bundle $E$ over a surface as above can be computed using the following proposition:
\begin{proposition}[See \cite{morita}, Proposition 4.11]\label{signature theorem}
Let $E$ be an oriented surface bundle over a surface, with closed oriented base and fibre. Then
$$3\sigma(E)=\langle e\cup e, [E]\rangle.$$
\end{proposition}

\section{Proof of Theorem \ref{thm:sign vol simp}}\label{preuve sign vol simp}

\begin{proof}[Proof of Theorem \ref{thm:sign vol simp}]
By Proposition \ref{signature theorem}, we have
$$\left\langle e\cup e, [E]\right\rangle=3\sigma(E).$$
On the other hand, by Proposition \ref{norme-vol simp},
$$|\left\langle e\cup e, [E]\right\rangle|=\|e\cup e\|\cdot\|E\|\leq \frac{1}{12}\|E\|,$$
as $\|e\cup e\|\leq \frac{1}{12}$ (see \cite{bucherConvPol}*{formula on p. 337}).
Hence $$\|E\|\geq 12\cdot 3 |\sigma(E)|=36|\sigma(E)|.$$
\end{proof}
In 1998, Kotschick proved the following theorem:
\begin{theorem}[\cite{kotschick}, Theorem 2]\label{kotschick-sigma chi}
Let $E$ be an aspherical surface bundle over a surface. Then
$$2|\sigma(E)|\leq \chi(E).$$
\end{theorem}
(Note that this is true even if $F$ or $B$ is the sphere, as the signature vanishes in these cases; see for example \cite{bryan-donagi-stipsicz}.)
The first author then obtained the following result:
\begin{theorem}[\cite{bucher}, Corollary 1.3 and \cite{bucherH^2xH^2}, Corollary 3]\label{michelle-vol simp fibres}
Let $F\hookrightarrow E \twoheadrightarrow B$ be an oriented surface bundle over a surface, with closed oriented base and fibre. Then

$$\|E\|\geq \|F\times B\|.$$

Furthermore, in the case of aspherical $F$ and $B$, the simplicial volume of the product $F\times B$ admits the value
$$\|F\times B\|=6\chi(F\times B).$$
\end{theorem}
Remember that 
$\chi(E)=\chi(F)\chi(B)=\chi(F\times B)$ for any $F$-bundle over $B$.
Putting everything together, we obtain:
$$\|E\|\geq \|F\times B\|=6\chi(E)\geq 12|\sigma(E)|.$$
In particular, the above inequality is weaker than the inequality of Theorem \ref{thm:sign vol simp}.

\section{Proof of Proposition \ref{prop: N et E}}\label{preuve N et E}
In this section we will work with group homology and cohomology, making use of the isometric isomorphisms
\[(H_*(E,\mathbb{R}), \|\cdot\|_1)\cong (H_*(\pi_1(E),\mathbb{R}), \|\cdot\|_1)\]
and
\[(H^*(E,\mathbb{R}), \|\cdot\|)\cong (H^*(\pi_1(E),\mathbb{R}), \|\cdot\|).\]
Slightly abusing notation, we will thus think of the fundamental class of $E$ as an element $[E]\in H_4(\pi_1(E), \mathbb{R})$ and of the Euler class of $E$ as a (bounded) element $e\in H^2(\pi_1(E), \mathbb{R})$.
Recall the alternation of a chain:
\begin{definition}
Let $\Gamma$ be a group, let $\underline{\gamma}=[\gamma_0, ..., \gamma_n ]\in C_n(\Gamma, \mathbb{R})$ be a basis element in the homogeneous chain complex of $\Gamma$. Define $\mathrm{Alt}(\underline{\gamma})$ by
$$\mathrm{Alt}(\underline{\gamma})=\frac{1}{(n+1)!}\sum_{\tau\in \mathrm{Sym}(n+1)}\signum(\tau)[\gamma_{\tau(0)}, ..., \gamma_{\tau(n)}]\in C_n(\Gamma, \mathbb{R}).$$
Denote by $\underline{\gamma}^\tau$ the element $[\gamma_{\tau(0)}, ..., \gamma_{\tau(n)}]$ obtained from $\underline{\gamma}$ by permuting the entries of $\underline{\gamma}$ by the permutation $\tau$. 

The definition is extended by linearity on the whole group $C_n(\Gamma, \mathbb{R})$.
\end{definition}
\begin{remark}\label{cycle et alterne}
It is well known that a cycle and its alternation define the same class in $H_n(\Gamma, \mathbb{R})$, that is $[z]=[\mathrm{Alt}(z)]\in H_n(\Gamma, \mathbb{R})$ (see \cite{fujiwara-manning}*{Appendix B} for a proof).
\end{remark}
\begin{remark}
Using the triangle inequality, one readily sees that for any $z\in C_n(\Gamma, \mathbb{R})$,
$$\|\mathrm{Alt}(z)\|_1\leq \|z\|_1.$$
\end{remark}

\begin{proof}[Proof of Proposition \ref{prop: N et E}]
Choose a fundamental cycle $\sum_{i=1}^k a_i\underline{\gamma^i}$ representing $[E]$. By definition, $[N]=e\cap[E]$.

Note that as $\|e\|\leq \frac{1}{2}$, we already have $2\|\left[N\right]\|_1\leq \|E\|$ by Proposition \ref{norme-vol simp}. With some more care, we will improve this by a factor $\frac{2}{3}$, getting the inequality of our proposition.

By Remark \ref{cycle et alterne} and the fact that the Euler class $e$ can be represented by $\frac{1}{2}\phi^*(Or)$, we have
$$[N]=e\cap [\mathrm{Alt}(E)]=\left[\sum_{i=1}^ka_i\frac{1}{5!}\sum_{\tau\in \mathrm{Sym}(5)}\signum(\tau)\frac{1}{2}\phi^*(Or)(\underline{\gamma^i}^\tau\rfloor)\lfloor \underline{\gamma^i}^\tau\right].$$

\begin{claim}\label{Claim} For every $i\in \{1, ..., k\}$ there exist chains $\alpha_i\in C_2(\pi_1(E))$ and $\beta_i\in C_3(\pi_1(E))$ such that
$$ \frac{1}{5!}\sum_{\tau\in \mathrm{Sym}(5)}\signum(\tau)\frac{1}{2}\phi^*(Or)(\underline{\gamma^i}^\tau\rfloor)\lfloor \underline{\gamma^i}^\tau = \alpha_i+\partial \beta_i,$$
with
$$\|\alpha_i\|\leq \frac{1}{3}.$$
 \end{claim}
 
Since $\sum_{i=1}^ka_i\alpha_i$ is by construction a cycle representing $[N]$, this immediately gives the inequality 
$$\| [N]\|_1 \leq \frac{1}{3} \sum_{i=1}^k|a_i|,$$
which, taking the infimum over all cycles representing $[E]$, proves the Proposition. It remains to prove Claim \ref{Claim}.

Fix $i\in \{1, ..., k\}$ and denote $\underline{\gamma^i}$ by $[\gamma_0, ..., \gamma_4]$.


Define 
$$T(j):= \frac{1}{2} \sum_{\tau\in \mathrm{Sym}(5),\\ \tau(2)=j}\signum(\tau)\frac{1}{2}\phi^*(Or)([\gamma_{\tau(0)}, \gamma_{\tau(1)}, \gamma_{\tau(2)}])\left[\gamma_{\tau(2)}, \gamma_{\tau(3)}, \gamma_{\tau(4)}\right]$$
for $j, \, 0\leq j\leq 4$, and note that $\frac{1}{5!}\sum_{\tau\in \mathrm{Sym}(5)}\signum(\tau)\frac{1}{2}\phi^*(Or)(\underline{\gamma^i}^\tau\rfloor)\lfloor \underline{\gamma^i}^\tau$ is equal to
$$\frac{1}{60}\left(T(0)+T(1)+T(2)+T(3)+T(4)\right).$$

Using that $e$ is alternating, we find that $2T(0)$ is a sum of at most $12$ signed singular simplices:
$$\begin{array}{l}
\phi^*(Or)([\gamma_3, \gamma_4, \gamma_0])\left[\gamma_0, \gamma_1, \gamma_2\right]-\phi^*(Or)([\gamma_2, \gamma_4, \gamma_0])\left[\gamma_0, \gamma_1, \gamma_3\right]\\
+\phi^*(Or)([\gamma_2, \gamma_3, \gamma_0])\left[\gamma_0, \gamma_1, \gamma_4\right]+\phi^*(Or)([\gamma_1, \gamma_2, \gamma_0])\left[\gamma_0, \gamma_3, \gamma_4\right]\\
-\phi^*(Or)([\gamma_1, \gamma_3, \gamma_0])\left[\gamma_0, \gamma_2, \gamma_4\right]+\phi^*(Or)([\gamma_1, \gamma_4, \gamma_0])\left[\gamma_0, \gamma_2, \gamma_3\right]\\
+\phi^*(Or)([\gamma_4, \gamma_3, \gamma_0])\left[\gamma_0, \gamma_2, \gamma_1\right]-\phi^*(Or)([\gamma_4, \gamma_2, \gamma_0])\left[\gamma_0, \gamma_3, \gamma_1\right]\\
+\phi^*(Or)([\gamma_3, \gamma_2, \gamma_0])\left[\gamma_0, \gamma_4, \gamma_1\right]+\phi^*(Or)([\gamma_2, \gamma_1, \gamma_0])\left[\gamma_0, \gamma_4, \gamma_3\right]\\
-\phi^*(Or)([\gamma_3, \gamma_1, \gamma_0])\left[\gamma_0, \gamma_4, \gamma_2\right]+\phi^*(Or)([\gamma_4, \gamma_1, \gamma_0])\left[\gamma_0, \gamma_3, \gamma_2\right].
\end{array}$$

Applying the cyclic permutation $(0, 1, 2, 3, 4)$ and its powers to the indices in the latter expression, we successively obtain $2T(1), 2T(2), 2T(3), 2T(4)$.



To prove the Claim, it is enough to prove it for $2T(0)$, that is to show that $2T(0)$ can be written, up to a boundary, as a sum of at most $\frac{2}{3}12=8$ singular simplices (with sign). 



Either all points $\phi(\gamma_0)x, ..., \phi(\gamma_4)x$ are distinct, or at least two of them coincide. We first consider the latter case.

Suppose two points coincide: if three or more points coincide, then at least $6$ terms vanish in $2T(0)$ and there is nothing to prove.
 
Suppose then that exactly two points coincide. Without loss of generality assume it is $\phi(\gamma_1)x$ with another one. If it coincides with $\phi(\gamma_0)x$, then $6$ terms in $2T(0)$ vanish. If not, we can by symmetry assume that it coincides with $\phi(\gamma_2)x$. By hypothesis, the points $\phi(\gamma_0)x,\phi(\gamma_1)x=\phi(\gamma_2)x, \phi(\gamma_3)x,\phi(\gamma_4)x$ are four distinct points on $S^1$. Up to reversing orientation we can suppose that $\phi(\gamma_0)x,\phi(\gamma_1)x=\phi(\gamma_2)x, \phi(\gamma_3)x$ are oriented positively and hence 
$$\phi^*(Or)([\gamma_1,\gamma_3,\gamma_0])=Or(\phi(\gamma_1),\phi(\gamma_3),\phi(\gamma_0))=+1.$$
\begin{figure}[htpb]
\centering
\begin{tikzpicture}
\draw(0,0) circle (1);
\fill (1,0) circle (2pt) node[right] {$\phi(\gamma_0)x$};
\fill(-0.7,0.7) circle (2pt) node[above left] {$\phi(\gamma_1)x=\phi(\gamma_2)x$};
\fill(-1,0) circle (2pt) node[left] {$\phi(\gamma_3)x$};
\end{tikzpicture}
\end{figure}

Using these explicit values for $\phi^*(Or)$ and the fact that the two terms where $\phi^*(Or)$ is evaluated on a triple containing $\gamma_1$ and $\gamma_2$ vanish in $2T(0)$, the expression for $2T(0)$ simplifies to
$$\begin{array}{l}
\phi^*(Or)([\gamma_3, \gamma_4, \gamma_0])\left[\gamma_0, \gamma_1, \gamma_2\right]-\phi^*(Or)([\gamma_2, \gamma_4, \gamma_0])\left[\gamma_0, \gamma_1, \gamma_3\right]+\left[\gamma_0, \gamma_1, \gamma_4\right]\\
-\left[\gamma_0, \gamma_2, \gamma_4\right]+\phi^*(Or)([\gamma_1, \gamma_4, \gamma_0])\left[\gamma_0, \gamma_2, \gamma_3\right]\\
+\phi^*(Or)([\gamma_4, \gamma_3, \gamma_0])\left[\gamma_0, \gamma_2, \gamma_1\right]-\phi^*(Or)([\gamma_4, \gamma_2, \gamma_0])\left[\gamma_0, \gamma_3, \gamma_1\right]-\left[\gamma_0, \gamma_4, \gamma_1\right]\\
+\left[\gamma_0, \gamma_4, \gamma_2\right]+\phi^*(Or)([\gamma_4, \gamma_1, \gamma_0])\left[\gamma_0, \gamma_3, \gamma_2\right].
\end{array}$$

Let us further examine the two possible values of $\phi^*(Or)([\gamma_3, \gamma_4, \gamma_0])$. If it is equal to $-1$, then we find in $2T(0)$ the expression
$$-\left[\gamma_0, \gamma_1, \gamma_2\right]+\left[\gamma_0, \gamma_1, \gamma_4\right]-\left[\gamma_0, \gamma_2, \gamma_4\right]=-\left[ \gamma_1,\gamma_2,\gamma_4\right]+\partial\left[ \gamma_0,\gamma_1,\gamma_2,v_4\right].$$

If $\phi^*(Or)([\gamma_3, \gamma_4, \gamma_0])=1$, then we can compute all values of $\phi^*(Or)$ since the vertices must be positioned as follows: 
\begin{figure}[htpb]
\centering
\begin{tikzpicture}
\draw(0,0) circle (1);
\fill (1,0) circle (2pt) node[right] {$\phi(\gamma_0)x$};
\fill(-0.7,0.7) circle (2pt) node[above left] {$\phi(\gamma_1)x=\phi(\gamma_2)x$};
\fill(-.433,-.933) circle (2pt) node[below left] {$\phi(\gamma_4)x$};
\fill(-1,0) circle (2pt) node[left] {$\phi(\gamma_3)x$};
\end{tikzpicture}
\end{figure}

The expression for $2T(0)$ thus becomes 
$$\begin{array}{l}
+\left[\gamma_0, \gamma_1, \gamma_2\right]-\left[\gamma_0, \gamma_1, \gamma_3\right]
+\left[\gamma_0, \gamma_1, \gamma_4\right]
-\left[\gamma_0, \gamma_2, \gamma_4\right]+\left[\gamma_0, \gamma_2, \gamma_3\right]\\
-\left[\gamma_0, \gamma_2, \gamma_1\right]+\left[\gamma_0, \gamma_3, \gamma_1\right]
-\left[\gamma_0, \gamma_4, \gamma_1\right]
+\left[\gamma_0, \gamma_4, \gamma_2\right]-\left[\gamma_0, \gamma_3, \gamma_2\right]
\end{array}$$
and we see that it contains 
$$\left[\gamma_0, \gamma_1, \gamma_2\right]-\left[\gamma_0, \gamma_1, \gamma_3\right]+\left[\gamma_0, \gamma_2, \gamma_3\right]=+\left[\gamma_1, \gamma_2, \gamma_3\right]-\partial \left[\gamma_0, \gamma_1,\gamma_2, \gamma_3\right].$$

Thus in both cases we have seen that the number of terms in $2T(0)$ can, up to taking boundaries, be reduced by at least a factor of $\frac{8}{12}=\frac{2}{3}$.

Suppose now the points $\phi(\gamma_0)x, ...,\phi(\gamma_4)x$ are all distinct. Without loss of generality, we suppose that $\phi(\gamma_0)x, ...,\phi(\gamma_4)x\in S^1$ are positively cyclically oriented.

\begin{figure}[htpb]
\centering
\begin{tikzpicture}
\draw(0,0) circle (1);
\fill (1,0) circle (2pt) node[right] {$\phi(\gamma_0)x$};
\fill(-0.7,0.7) circle (2pt) node[above left] {$\phi(\gamma_2)x$};
\fill(-.433,-.933) circle (2pt) node[below left] {$\phi(\gamma_4)x$};
\fill(0.433,0.933) circle (2pt) node[above right] {$\phi(\gamma_1)x$};
\fill(-1,0) circle (2pt) node[left] {$\phi(\gamma_3)x$};
\end{tikzpicture}
\end{figure}

The expression $2T(0)$ simply becomes
$$\begin{array}{l}
+\left[\gamma_0, \gamma_1, \gamma_2\right]-\left[\gamma_0, \gamma_1, \gamma_3\right]
+\left[\gamma_0, \gamma_1, \gamma_4\right]+\left[\gamma_0, \gamma_3, \gamma_4\right]\\
-\left[\gamma_0, \gamma_2, \gamma_4\right]+\left[\gamma_0, \gamma_2, \gamma_3\right]
-\left[\gamma_0, \gamma_2, \gamma_1\right]+\left[\gamma_0, \gamma_3, \gamma_1\right]\\
-\left[\gamma_0, \gamma_4, \gamma_1\right]-\left[\gamma_0, \gamma_4, \gamma_3\right]
+\left[\gamma_0, \gamma_4, \gamma_2\right]-\left[\gamma_0, \gamma_3, \gamma_2\right]
\end{array}$$
and we use the two cycle relations 
$$\begin{array}{l}
+\left[\gamma_0, \gamma_1, \gamma_2\right]-\left[\gamma_0, \gamma_1, \gamma_3\right]+\left[\gamma_0, \gamma_2, \gamma_3\right]=+\left[\gamma_1, \gamma_2, \gamma_3\right]-\partial \left[\gamma_0, \gamma_1,\gamma_2, \gamma_3\right]\\
-\left[\gamma_0, \gamma_2, \gamma_1\right]+\left[\gamma_0, \gamma_3, \gamma_1\right]-\left[\gamma_0, \gamma_3, \gamma_2\right]=-\left[ \gamma_3, \gamma_2,\gamma_1\right]+\partial\left[\gamma_0, \gamma_3, \gamma_2,\gamma_1\right]
\end{array}$$
to conclude as before that the number of terms can up to taking boundaries be reduced by a factor of $\frac{2}{3}$, which finishes the proof of the Claim and the Proposition.\end{proof}

\begin{remark}
Note that we have $\pi_*(\left[N\right])=\chi(F)\left[B\right]$, and consequently $$2\chi(E)=|\chi(F)|\|B\|\leq\|\left[N\right]\|_1.$$ Hence for bundles with $\|E\|=6\chi(E)$, we obtain the equality $$\|\left[N\right]\|_1=2\chi(E)=\frac{1}{3}\|E\|.$$
This includes all bundles with finite image of the holonomy homomorphism, and in particular the trivial bundle $E=F\times B$, for which $[N]$ can be represented by $\chi(F)$ disjoint copies of $B$ in $F\times B$.
\end{remark}
\section{Ramified coverings}\label{ramifie}
In this section we present a method for constructing surface bundles with non-zero signature using ramified coverings and then study the simplicial volume of the total space of such bundles.
\subsection{Construction of surface bundles using ramified coverings}\label{rev ramifies}
The first examples of surface bundles over surfaces with non-zero signature were constructed  independently by Kodaira \cite{kodaira} in 1967 and Atiyah \cite{atiyah} in 1969 with a method relying on ramified coverings. We outline this method here, following its exposition in \cite{morita}*{Paragraph 4.3.3}.

First choose a closed oriented surface $\Sigma^o=\Sigma_{g_0}$, with genus $g_0\geq 2$. Then take a $d$-fold cyclic covering $\rho\colon \Sigma\rightarrow \Sigma^o$ of $\Sigma^o$, and let $\sigma$ be a generator of its covering transformation group $\mathbb{Z}/d\mathbb{Z}$.
\begin{remark}
This implies that $\sigma^i$ is fixed point free, for $1\leq i\leq d-1$.
\end{remark}
Denote by $g$ the genus of $\Sigma$. We have $2-2g=d(2-2g_0)$.
Consider the following homomorphisms:
$$\pi_1(\Sigma)\longrightarrow \pi_1(\Sigma)^{ab}\cong H_1(\Sigma, \mathbb{Z})\longrightarrow H_1(\Sigma, \mathbb{Z}/d\mathbb{Z})\cong (\mathbb{Z}/d\mathbb{Z})^{2g}.$$
Their composition is surjective and its kernel is a normal subgroup of finite index in $\pi_1(\Sigma)$. As such, it defines a finite regular covering $\rho'\colon \Sigma'\rightarrow \Sigma$.
We have a map $\sigma^i\circ\rho'\colon \Sigma'\rightarrow \Sigma$ for each $1\leq i\leq d$. We can then consider the graph of $\sigma^i\circ\rho'$ in $\Sigma'\times \Sigma$ for each $i$: it defines a submanifold $\Gamma_{\sigma^i\circ\rho'}$.
\begin{remark}\label{rem: disjoint union}
The fact that $\sigma^i$ is fixed point free for all $1\leq i\leq d-1$ ensures that the graphs $\Gamma_{\sigma^i\circ\rho'}, \Gamma_{\sigma^j\circ\rho'}$ are disjoint whenever $i\neq j$.
\end{remark}
Take the disjoint union $\Gamma_{\sigma\circ\rho'}\sqcup...\sqcup\Gamma_{\sigma^d\circ\rho'}$ of these submanifolds and denote it by $D$. It is of codimension $2$ in $\Sigma'\times \Sigma$, therefore it defines a class $[D]\in H_2(\Sigma'\times \Sigma, \mathbb{Z})$.

We will need the following proposition:
\begin{proposition}[\cite{morita}, Proposition 4.10]\label{classe divisible}
Let $B$ be a closed oriented $C^\infty$ manifold and let $D\subset B$ be an oriented submanifold of codimension $2$. Suppose that, for some $m\in\mathbb{Z}_{> 0}$, the homology class $[D]\in H_{n-2}(B, \mathbb{Z})$ determined by $D$ is divisible by $m$ in $H_{n-2}(B, \mathbb{Z})$. Then there exists an $m$-fold cyclic ramified covering $\widetilde{B}\rightarrow B$ ramified along $D$.
\end{proposition}
The class $[D]$ defined above is divisible by $d$ in $H_2(\Sigma'\times \Sigma, \mathbb{Z})$ \cite{morita}*{p. 158}. Thus using Proposition \ref{classe divisible}, we obtain a ramified covering $f\colon E\rightarrow \Sigma'\times \Sigma$ of degree $d$ ramified along $D$.

Finally we get a surface bundle $E\rightarrow \Sigma'$ as the composition $E\stackrel{f}{\rightarrow} \Sigma'\times \Sigma\rightarrow \Sigma'$, where $\Sigma'\times \Sigma\rightarrow \Sigma'$ is the canonical projection to the first factor. The fibre of $E$ is $f^{-1}(\Sigma)$.

The signature of $E$ can be explicitly computed and it is non-zero.
For this, one more result is used, giving relations between the Euler class of $E$ and that of $\Sigma'\times \Sigma$.
\begin{proposition}[\cite{morita}, Proposition 4.12]\label{f*nu, e tilde}
Let $\pi\colon E\rightarrow B$ and $\widetilde{\pi}\colon \widetilde{E}\rightarrow B$ be two surface bundles over the same base space $B$. Suppose that there is a map $f\colon \widetilde{E}\rightarrow E$ between the total spaces which is a $d$-fold cyclic ramified covering ramified along an oriented submanifold $D\subset E$ of codimension $2$, and that $f$ is a bundle map (i. e. $\pi\circ f=\widetilde{\pi}$). Suppose also that $D$ intersects each fibre of $\pi$ transversely at exactly $d$ points, and write $\widetilde{D}=f^{-1}(D)$.
$$
\xymatrix
{
\widetilde{D}\ar@{^{(}->}[d]& D\ar@{^{(}->}[d]\\
\widetilde{E}\ar[r]^f\ar[d]^{\widetilde{\pi}}& E\ar[d]^\pi\\
B\ar[r]^=&B
}
$$Then:
\begin{enumerate}
\item $f^*(\nu)=d\tilde{\nu}$;
\item $\tilde{e}=f^*\left(e-(1-\frac{1}{d})\nu\right)$\label{euler class ramified},
\end{enumerate}
where $\nu$, respectively $\tilde{\nu}$, represents the Poincar\'e dual of the homology class of $D$, respectively $\widetilde{D}$, and $e$, respectively $\tilde{e}$, denotes the Euler class of $\pi$, respectively $\tilde{\pi}$.
\end{proposition}
All the assumptions of Proposition \ref{f*nu, e tilde} are satisfied by $f\colon E\rightarrow \Sigma'\times \Sigma$.

\subsection{Simplicial volume of such bundles}\label{ramifie calcul}

As observed at the end of Section \ref{preuve sign vol simp}, the results of the first author show that 
$$\|E\|\geq 6\chi(E)$$
for any aspherical surface bundle $E$ over a surface \cites{bucherH^2xH^2, bucher}.

Now if we restrict our attention to surface bundles over surfaces coming from the ramified covering construction explained in the previous subsection, we can enhance this inequality. 

Consider the following diagram that represents the aforementioned construction:
$$\xymatrix
{
E\ar[r]^-f&\Sigma'\times\Sigma\ar[dr]_{p}\ar[dl]^{p'}\\
\Sigma' &&\Sigma
}
$$
The maps $p$ and $p'$ are the natural projections.
The map $f$ is a cyclically ramified covering of degree $d$ of $\Sigma'\times\Sigma$, ramified along the codimension $2$ submanifold $D\subset \Sigma'\times \Sigma$ defined above, and $\Sigma'$ is a $d'$-fold covering of $\Sigma$. The intersection (both algebraic and geometric) $D\cap \Sigma'$ in $\Sigma'\times \Sigma$ consists of $d'd$ points while the intersection $D\cap \Sigma$ consists of $d$ points.
\begin{remark}
In order to avoid heavy notation, by $\Sigma\subset \Sigma'\times \Sigma$ we mean the choice of a subsurface $\{x'\}\times\Sigma$. Similarly $[\Sigma]\in H^2(\Sigma'\times \Sigma, \mathbb{Z})$ denotes a class $[\{x'\}\times\Sigma]$.

For further use, we also mark that the notation $[A]^*$ stands for the Poincar\'e dual of the homology class $[A]$.
\end{remark}
The crucial remark, already made by Bryan, Donagi and Stipsicz in \cite{bryan-donagi-stipsicz} and LeBrun in \cite{lebrun}, is that $E$ admits (at least) two different bundle structures: namely the compositions $p\circ f$ and $p'\circ f$ are the bundle projections of the surface bundles $\pi\colon E\rightarrow\Sigma$ and $\pi'\colon E\rightarrow \Sigma'$ with fibres $f^{-1}(\Sigma')$ and $f^{-1}(\Sigma)$ respectively.
\newtheorem*{thm:fibre ramifie}{Theorem \ref{thm:fibre ramifie}}

\begin{proof}[Proof of Theorem \ref{thm:fibre ramifie}]
Denote by $e=\chi(\Sigma')[\Sigma]^*$ the Euler class of the product bundle $\Sigma'\times\Sigma\rightarrow\Sigma$, and by $e'=\chi(\Sigma)[\Sigma']^*$ the Euler class of the product bundle $\Sigma'\times\Sigma\rightarrow\Sigma'$, both in $H^2(\Sigma'\times\Sigma, \mathbb{Z})$. 

By Proposition \ref{f*nu, e tilde}\ref{euler class ramified}, the Euler class of the bundle $\pi$ is $e_E=f^*\left(e-(1-\frac{1}{d})[D]^*\right)$ and the one of the bundle $\pi'$ is $e_E'=f^*\left(e'-(1-\frac{1}{d})[D]^*\right)$, both in $H^2(E, \mathbb{Z})$.

We compute:
\begin{eqnarray*}
\left\langle e_E'\cup e_E, [E]\right\rangle&=&\left\langle f^*\left(e'-(1-\frac{1}{d})[D]^*\right)\cup f^*\left(e-(1-\frac{1}{d})[D]^*\right), [E]\right\rangle\\
						&=&d\left\langle \left(e'-(1-\frac{1}{d})[D]^*\right)\cup \left(e-(1-\frac{1}{d})[D]^*\right), [\Sigma'\times\Sigma]\right\rangle\\
						&=&d\left\langle e'-(1-\frac{1}{d})[D]^*, \chi(\Sigma')[\Sigma]-(1-\frac{1}{d})[D]\right\rangle\\
						&=&d\left(\chi(\Sigma')\chi(\Sigma)-\chi(\Sigma')(1-\frac{1}{d})d\right. \\
						& &\left. -(1-\frac{1}{d})\chi(\Sigma)dd'+(1-\frac{1}{d})^2dd'\chi(\Sigma)\right)\\
						&=&d\left(\chi(\Sigma')\chi(\Sigma)-\chi(\Sigma')(d-1)-\chi(\Sigma')(d-1)+(1-\frac{1}{d})^2d\chi(\Sigma')\right)\nonumber\\
						&=&d\left(\chi(\Sigma')\chi(\Sigma)-2\chi(\Sigma')(d-1)+(1-\frac{1}{d})^2d\chi(\Sigma')\right)\\
						&=&d\chi(\Sigma')\left(\chi(\Sigma)-2(d-1)+(d-1)(1-\frac{1}{d})\right)\\
						&=&d\chi(\Sigma')\left(\chi(\Sigma)-(d-1)(1+\frac{1}{d})\right).
\end{eqnarray*}
Note that we used $[D]^*\cap[D]=dd'\chi(\Sigma)$, which can be proven as follows: Let us denote by $\Delta$ the diagonal in $\Sigma\times\Sigma$, that is the set $\left\{(x, x)\in\Sigma\times \Sigma\,\mid\,x\in\Sigma\right\}$. It defines a class $[\Delta]\in H_2(\Sigma\times\Sigma,\mathbb{Z})$. Using Remark \ref{rem: disjoint union}, we compute
\begin{eqnarray*}
[D]^*\cap[D]&=&\left([\Gamma_{\sigma\circ\rho'}]+...+[\Gamma_{\sigma^d\circ\rho'}]\right)^*\cap\left([\Gamma_{\sigma\circ\rho'}]+...+[\Gamma_{\sigma^d\circ\rho'}]\right)\\
	&=&\sum_{i=1}^d[\Gamma_{\sigma^i\circ\rho'}]^*\cap[\Gamma_{\sigma^i\circ\rho'}]\\
	&=&d\,[\Gamma_{\rho'}]^*\cap[\Gamma_{\rho'}]\\
	&=&dd'\,[\Delta]^*\cap[\Delta]\\
	&=&dd'\chi(\Sigma).
\end{eqnarray*}
This computation can also be extracted from \cite{morita}*{pp. 157-160}. The last equality is a general fact about the diagonal in a product of manifolds, see for example \cite{bott-tu}*{p. 128}.

By \cite{bucher}*{Proposition 2.1} and because $\|e_E\|\leq \frac{1}{2}$ by Subsection \ref{def classe euler}, we have:
$$\|e_E'\cup e_E\|\leq \frac{1}{3}\|e_E\|\leq \frac{1}{6}.$$
So we obtain:
\begin{eqnarray}
|\langle e_E'\cup e_E, [E]\rangle|&=&d|\chi(\Sigma')||\chi(\Sigma)-(d-1)(1+\frac{1}{d})|,\nonumber\\
\frac{1}{6}\|E\|\geq \|e_E'\cup e_E\|\|E\|&=&d|\chi(\Sigma')|\left(|\chi(\Sigma)|+(d-1)(1+\frac{1}{d})\right),\nonumber\\
\|E\|&\geq& 6d|\chi(\Sigma')|\left(|\chi(\Sigma)|+(d-1)(1+\frac{1}{d})\right).\nonumber
\end{eqnarray}
The fibre $f^{-1}(\Sigma)$ of the bundle $\pi'$ has Euler characteristic
$$\chi(f^{-1}(\Sigma))=d\chi(\Sigma)-d(d-1),$$
as it is a degree $d$ cyclic ramified covering of $\Sigma$ with $d$ intersection points with the ramification locus $D$. The Euler characteristic of $E$ can then be written as
$$\chi(E)=\chi(\Sigma')\left(d\chi(\Sigma)-d(d-1)\right)=d|\chi(\Sigma')|\left(|\chi(\Sigma)|+(d-1)\right).$$
The result can thus be expressed as
$$\|E\|\geq 6\chi(E)+6|\chi(\Sigma')|(d-1).$$
\end{proof}





\end{document}